\documentclass{article}
\usepackage{graphicx}
\usepackage{latexsym}
\usepackage{amsmath}
\usepackage{amsthm}
\usepackage{subfig}
\usepackage{url}
\usepackage{bm}

\usepackage{amssymb}
 
\newcommand{\vn}{{\bf n}}

\newcommand{\vQ}{{\bf Q}}

\date{\today}
\newcommand{\lbl}[1]{\label{#1}}

\newcommand{\be}{\begin{eqnarray}}
\newcommand{\ee}{\end{eqnarray}}

\newcommand{\R}{\mathbb{R}}
\newcommand{\half}{\frac{1}{2}}

\newcommand{\om}{\Omega}

\newcommand{\1}{{\bf 1}}

\newcommand{\n}{{\bf n}}
\newcommand{\x}{{\bf x}}
\newcommand{\m}{{\bf m}}
\newcommand{\Q}{{\bf Q}}
\newcommand{\Pz}{{\bf P}}
\newcommand{\az}{{\bf a}}

\newcommand{\e}{{\bf e}}
\newcommand{\uu}{{\bf u}}
\newcommand{\y}{{\bf y}}
\newcommand{\A}{{\bf A}}
\newcommand{\D}{{\bf D}}
\newcommand{\Rz}{{\bf R}}
\newcommand{\nnu}{{\bm \nu}}
\newtheorem{thm}{Theorem}

\numberwithin{equation}{section}

\newcommand{\vertiii}[1]{{\bm \|} #1 
    {\bm \|}}
\graphicspath{%
    {converted_graphics/}
    {C:/Users/John/Documents/Papers/}
    {/}
}
\begin{document}
\title
{Discontinuous order parameters in liquid crystal theories}
\author{J. M. Ball and S. J. Bedford\\ \\
Mathematical Institute,   University of Oxford,\\
Andrew Wiles Building,
Radcliffe Observatory Quarter,\\
Woodstock Road,
Oxford,
OX2 6GG,
 U.K. }

 \maketitle
 \markboth{ }{}
\begin{abstract}   
The paper is concerned with various issues surrounding the mathematical description of defects in models of liquid crystals, drawing on experience from solid mechanics. The roles played by a suitable choice of function space and by the growth properties of the free-energy density are highlighted. Models in which the director can jump across surfaces are formulated, and their relevance for nematic elastomers, order reconstruction and smectic A thin films discussed.
\end{abstract}

\section{Introduction}
\setcounter{equation}{0}
In this paper we consider some physical and mathematical issues involved in addressing the following questions:
\begin{itemize}\item Should we allow continuum order parameters (directors, $\Q$-tensors etc) describing liquid
crystals to jump across surfaces?
\item Is there a useful theory of nematic and cholesteric
 liquid crystals in which the order parameter is a unit vector field
and for which observed defects have finite energy?
\item How should we choose an appropriate function space for mathematical models of liquid crystals? 
\end{itemize}
In exploring these questions we draw on lessons from the study of related questions from solid mechanics, in which the central model of nonlinear elasticity has a similar variational structure to models of liquid crystals.  In solid mechanics defects involving discontinuities across surfaces (fracture surfaces and phase boundaries) are commonplace, and their study leads to similar issues regarding problem formulation, function spaces and growth conditions of energy densities. Surface discontinuities are not usually considered to be relevant for liquid crystals. Nevertheless there are situations in which such discontinuities appear to arise, in liquid crystal elastomers, order reconstruction problems and smectic thin films, and they can also be useful in a purely mathematical context as a device for handling nonorientable director fields. These examples are treated with varying degrees of detail below.

Although we have tried to make the statements in the paper mathematically accurate, we have chosen to omit the often rather technical mathematical proofs of theorems, these being more suitable for a specialist mathematical readership (see \cite{u14,bedfordcholesterics, bedfordthesis}). 

\section{Classical models of liquid crystals}
To set the scene, we recall three classical models of static configurations of nematic and cholesteric liquid crystals occupying a bounded open region $\om\subset\R^3$ with sufficiently smooth boundary $\partial\Omega$\footnote{For the validity of the theorems stated in the paper it is enough to suppose that $\Omega$ is a Lipschitz domain, examples of which include balls and cubes.}   described in terms of their corresponding free-energy functionals. The oldest of these is the {\it Oseen-Frank theory}, in which the order parameter is a unit vector director field $\n=\n(\x)$ with corresponding free-energy functional
\be
\lbl{oseenfrank}
I_{OF}(\n)=\int_\om W(\n,\nabla \n)\,d\x,
\ee
where 
\be 
\lbl{ofw}
W(\n,\nabla \n)=K_{1}(\textrm{div}\,
\n)^{2} + K_{2}(\n\cdot \textrm{curl}\,\n+q_0)^{2}  +
K_{3}|\n \times \textrm{curl}\,\n|^{2}\nonumber\\&& \mbox{ }\hspace{-3in} +
(K_{2}+K_{4})(\textrm{tr}(\nabla \n)^{2} - (\textrm{div}
\n)^{2}),
\ee
and the coefficients $K_i, q_0$ are constants, with $q_0=0$ for nematics and $q_0\neq 0$ for cholesterics. (The constants $K_i, q_0$, together with other material constants specified below,  depend on temperature, but for simplicity we do not make this explicit and assume we are working at a constant temperature at which the liquid crystal is in the nematic phase.) 

In the {\it Landau - de Gennes theory} the order parameter is a symmetric traceless $3\times 3$ matrix $\Q=\Q(\x)=(Q_{ij}(\x))$ with corresponding free-energy functional
\be 
\label{ldg}
I_{LdG}(\Q)=\int_\om \psi(\Q,\nabla \Q)\,d\x.
\ee
The free-energy density can be written as $\psi(\Q,\nabla\Q)=\psi_B(\Q)+\psi_E(\Q,\nabla\Q)$, where $\psi_B(\Q)=\psi(\Q,0)$. The bulk energy $\psi_B$ is by frame-indifference an isotropic function of $\Q$, so that 
\be\label{isotropy}
\psi_B(\Q)=\hat\psi_B({\rm tr}\,\Q^2, \det\Q)
\ee
for some function $\hat\psi_B$, and is
often assumed to have the quartic form
\be 
\label{bulk}\psi_B(\Q)=a\,{\rm tr}\,\Q^2-\frac{2b}{3}{\rm tr}\,\Q^3+c\,{\rm tr}\,\Q^4,
\ee
where $b>0,c>0$ and  $a$ are constants. A possible frame-indifferent form for the elastic energy is 
\be 
\label{elastic}
&&\psi_E(\Q,\nabla \Q)=\sum_{i=1}^5L_iI_i,
     \\ &&I_1 =   Q_{ij,j}  Q_{ik,k}, \;
    I_2 =   Q_{ik,j}Q_{ij,k},\nonumber\\ 
&&I_3 =  Q_{ij,k}Q_{ij,k}, \;
    I_4 = Q_{lk} Q_{ij,l} Q_{ij,k},\; I_5=\varepsilon_{ijk}Q_{il}Q_{jl,k},\nonumber
\ee
where the $L_i$ are  material parameters. The bulk energy \eqref{bulk} has the drawback that there is no mechanism in the resulting theory for ensuring   the preservation of the constraint that the minimum eigenvalue $\lambda_{\rm min}(\Q)$ of $\Q$ is greater than $-\frac{1}{3}$, which arises from the interpretation of the $\Q$-tensor as a normalized probability distribution of molecular orientations. The natural way to preserve this constraint is via a bulk energy satisfying 
\be 
\label{blowup}
\psi_B(\Q)\to\infty\mbox{  as  }\lambda_{\rm min}\to -\frac{1}{3},
\ee
 as was proposed by Ericksen \cite{ericksen1991liquid} in the context of his theory described below. Such a bulk-energy function has been derived on the basis of the Onsager theory with the Maier-Saupe potential in 
\cite{katrieletal,u9,j59}. 

The Oseen-Frank energy \eqref{oseenfrank} can be obtained formally  from the Landau - de Gennes energy \eqref{ldg} (with the addition of the constant $q_0^2$) by making the ansatz that $\Q$ is uniaxial with constant scalar order parameter $s$, namely
\be 
\label{uniaxial}
\Q(\x)=s\left(\n\otimes\n-\frac{1}{3}{\bf 1}\right),
\ee
motivated by the fact that the bulk energy $\psi_B(\Q)$ given by \eqref{bulk} is minimized by $\Q$ of the form \eqref{uniaxial} with $s=0$ (the isotropic state) if $a\geq \frac{b^2}{27c}$ and   $s=\frac{b+\sqrt{b^2-24ac}}{4c}$ (the nematic state) when $a\leq \frac{b^2}{27c}$. If $a=\alpha(\theta-\theta^*)$ where $\theta$ is the temperature and $\alpha>0$, as is commonly assumed, then the condition  $a\leq \frac{b^2}{27c}$ reduces to $\theta\leq\theta_c:=\theta^*+\frac{b^2}{27\alpha c}$. In this case the $K_i$ are given by the explicit formulae \cite{morietal, mottramnewton}
\begin{eqnarray}  
&K_1=(L_1+L_2+2L_3)s^2-\frac{2}{3}L_4s^3,& K_2=2L_3s^2-\frac{2}{3}L_4s^3,\nonumber\\&K_3=(L_1+L_2+2L_3)s^2+\frac{4}{3}L_4s^3,&
K_4=L_2s^2,\\
&\hspace{-1in}q_0=\frac{L_5}{4(L_3+\frac{2}{3}L_4)}.\nonumber
\end{eqnarray} 
 Thus   $I_{OF}(\n)$ may be regarded as the energy corresponding to $I_{LdG}(\Q)$ when the elastic constants $L_i$ are small, as analyzed rigorously by \cite{majumdarzarnescu}, \cite{nguyenzarnescu} for the one-constant nematic case   $L_1=L_2=L_4=0,\;L_3>0,\; q_0=0$ with $\psi_B$ given by \eqref{bulk}. 

In the {\it Ericksen theory}   \cite{ericksen1991liquid} the more general ansatz is assumed in which $\Q$ has the uniaxial form \eqref{uniaxial} but   $s$ is allowed to depend on $\x$, leading to a theory in which the  order parameter is the pair $(s, \n)$ with corresponding free-energy functional
\be 
\label{ericksen}
I_E(s,\n)=\int_\om W(s,\nabla s, \n, \nabla\n)\,d\x.
\ee
Thus in the Ericksen theory the 5-dimensional order parameter of the Landau - de Gennes theory is reduced to a 3-dimensional order parameter, while in the Oseen-Frank theory it is reduced to a 2-dimensional order parameter. 

\section{Function spaces}

It is an uncomfortable fact that, as well as giving the free-energy functional, we  need to specify a {\it function space} to which the order parameter belongs, that expresses the worst kind of singularities that the order parameter is allowed to have.  This leads to some inescapable technical mathematical issues that are necessary so as to  define precisely the model being studied.  For example, in the case of the Oseen-Frank theory it would not make sense to suppose that $\n$ is a continuously differentiable unit vector field in $\om$, because we know that the Oseen-Frank theory supports singular director configurations, such as the hedgehog
\begin{equation}\label{hedgehog}
\tilde\n(\x)=\frac{\x}{|\x|},
\end{equation}
which are not continuously differentiable and which correspond to observed defects. Expressed in terms of function spaces, $\n$ given by \eqref{hedgehog} does not belong to $C^1(\om;{\mathbb R}^3)$, where for $r=1,2,\ldots$,
\begin{eqnarray*}C^r(\Omega;{\mathbb R}^m)= \{ r \mbox{ times continuously} \mbox{ differentiable maps }\uu:\Omega\rightarrow{\mathbb R}^m\}.\end{eqnarray*}
The usual function spaces that are used to handle variational problems such as \eqref{oseenfrank}, \eqref{ldg}, \eqref{ericksen} are the Sobolev spaces defined for $1\leq p\leq\infty$ by
\begin{eqnarray*}
W^{1,p}(\Omega;{\mathbb R}^m)=\{\uu\in L^p(\Omega;{\mathbb R}^m):
 \nabla \uu\in L^p(\Omega;M^{m\times n})\}.
\end{eqnarray*}
Here $\om\subset\R^n$ is a bounded domain (in this paper $n=3$ always), $L^p(\Omega;{\mathbb R}^m)$ denotes the space of (measurable) mappings $\uu:\om\to{\mathbb R}^m$ such that $\|\uu\|_p<\infty$, where
\begin{eqnarray*}\|\uu\|_p=\left\{\begin{array}{ll}\left(\int_\Omega |\uu|^p\,d\x\right)^{\frac{1}{p}}&\mbox{if }1\leq p<\infty\\
\mbox{ess sup}_{\x\in\Omega}|\uu(\x)|&\mbox{if }p=\infty
\end{array}\right.,
\end{eqnarray*}
$M^{m\times n}$ denotes the set of real $m\times n$ matrices, and $\mbox{ess sup}$ means the supremum disregarding subsets of $\om$ of zero ($n$-dimensional) volume.
The derivative $\nabla \uu(\x)=\left(\frac{\partial u_i}{\partial x_\alpha}(\x)\right)$, which belongs to $M^{m\times n}$  for each $\x$, is the
  weak or distributional derivative of $\uu$, which is defined so that the formula 
\begin{eqnarray*}
\int_\Omega \frac{\partial u_i}{\partial x_\alpha}\varphi\,d\x=-\int_\Omega u_i\frac{\partial\varphi}{\partial x_\alpha}\,d\x 
\end{eqnarray*}
for integration by parts holds for any smooth test function $\varphi:\om\to{\mathbb R}^m$ that vanishes together with all its derivatives in a neighbourhood of the boundary $\partial\om$ of $\om$.  In the applications to liquid crystals the mapping $\uu$ may be  $\n$, $\Q$ or the pair $(s,\n)$ depending on the model used. In the case of the director $\n$, which takes values in the unit sphere $S^2$, we will use the shorthand $W^{1,p}(\om;S^2)=\{\n\in W^{1,p}(\om;{\mathbb R}^3): |\n(\x)|=1\mbox{ for almost every } \x\in\om\}$, where `almost every' means for every $\x$ except possibly for a set of points of zero   volume. A somewhat fuller introduction to Sobolev spaces for those working in liquid crystals is given in \cite{j57}, while for a more mathematical description the reader may, for example, consult  \cite{evanspde}.

 Mappings $\uu$ that are discontinuous across a  two-dimensional surface (for example, a plane) do not belong to any of the Sobolev spaces. In order to treat such mappings it has become standard to use spaces of mappings of bounded variation, which admit mappings that are discontinuous across a surface while at the same time allowing them to have a well-defined derivative on each side of the  surface. The particular space of this type that we will use in this paper is the space $SBV(\om;{\mathbb R}^m)$ of special mappings of bounded variation. This consists of those mappings $\uu:\om\to{\mathbb R}^m$ in $L^1(\om;{\mathbb R}^m)$, that is
\begin{equation}\label{2.6}
 \int_\Omega |{\bf u}(\x)|\,d\x<\infty,
\end{equation}
whose (distributional) derivative   is a {\it measure} having no Cantor part. A measure has a more broad definition than a normal mapping, allowing, for example, for Dirac masses,  and is well described in \cite{evans92}. The exclusion of those $\uu$ whose derivative  has a Cantor part  is to eliminate such pathological examples as the Cantor function \cite{cantor} or `devil's staircase', which is a continuous nondecreasing function $f:[0,1]\to [0,1]$ whose  derivative $f'$ is zero at almost every point of $[0,1]$, but is such that $f(0)=0$ and $f(1)=1$, so that the fundamental theorem of calculus $\int_0^1f'(x)\,dx=f(1)-f(0)$ fails. The precise definition of $SBV(\om;{\mathbb R}^m)$ can be found in \cite{ambrosioetal00}. If $\uu\in SBV(\om;{\mathbb R}^m)$ then  almost every point of $\om$ is a point of `approximate differentiability' of $\uu$, at which an approximate gradient $\nabla\uu$ is defined, while there is a {\it jump set} $S_\uu$ where $\uu$ has approximate jump discontinuities and which has a special structure, with a well-defined normal and approximate one-sided limits of $\uu$ except for a set of points of zero two-dimensional Hausdorff measure (area). We shall use the notation $SBV(\om;S^2)$ to describe the space of $\n\in SBV(\om;\R^3)$ with $|\n(\x)|=1$ almost everywhere. We also use the notation $SBV_{\rm loc}(\Omega;S^2)$ to denote the space of $\n$ belonging to $SBV(U;S^2)$ for any open subset $U\subset\om$ which is a positive distance from the boundary $\partial\om$.

\section{Lessons from solid mechanics}
\label{solids}
Nonlinear elasticity has a similar variational structure to continuum
models of liquid crystals, with a free-energy functional of the form
\be \label{elasticity}I(\y)=\int_\om \psi(\nabla\y(\x))\,d\x,\ee
where $\y(\x)\in\R^3$ is the deformed position of a material point having 
position $\x$ in a reference configuration $\om\subset\R^3$. Minimizers of \eqref{elasticity} can have singularities, and the predictions of the model depend on the function space used.

As an example, consider the problem of minimizing  for given $\lambda>0$ the free energy
\be \label{rubber}I(\y)=\int_{B(0,1)}(|\nabla\y|^2+h(\det\nabla\y))\,d\x\ee
among   $\y$ satisfying the boundary condition $\y(\x)=\lambda \x$ for  $|\x|=1$, where $B(0,1)=\{\x\in {\mathbb R}^3:|\x|<1\}$ and $h:\R\to [0,\infty]$ satisfies
\begin{eqnarray*}
 h(\delta)=\infty \mbox{ for }   \delta\leq 0, \;\;
h \mbox{ finite, smooth and strictly convex for }\delta>0,\\ h(\delta)\to \infty\mbox{ as }\delta\to 0+,\;\;
\frac{h(\delta)}{\delta}\to\infty \mbox{ as } \delta\to \infty.
\end{eqnarray*}
(Here  and below we use the notation $|\A|=({\rm tr}\,\A^T\A)^\frac{1}{2}$ for the Euclidean norm of a $3\times 3$ matrix $\A$.)
This functional is a compressible form of the widely used neo-Hookean model\footnote{The integrand in \eqref{rubber} is {\it polyconvex}, that is it is a convex function of the minors  (subdeterminants) of $\nabla\y$. Polyconvexity is a sufficient condition for $\psi$ to be {\it quasiconvex}, a property that is known to be the central convexity condition of the multi-dimensional calculus of variations but which is poorly understood. The free-energy functions usually used for liquid crystals are convex, whereas  one would expect on mathematical grounds that quasiconvex free-energy functions would arise. For a survey of these different convexity condtions see \cite{dacorogna}.} of rubber. The condition that $h(\delta)=\infty$ for $\delta\leq 0$ enforces the non-interpenetration constraint $\det \nabla\y(\x)>0$ for almost every $\x$.  The unique minimizer among smooth mappings, or in $W^{1,p}(\om;\R^3)$ for $p\geq 3$ is the homogeneous deformation $\y^*(\x)=\lambda \x$ (this is because $\psi$ is quasiconvex --- see the preceding footnote). However $\y^*$ is not the minimizer for sufficiently large $\lambda$ in the natural energy space $W^{1,2}(\om;\R^3)$, since the material can reduce its energy by creating one or more cavities. In fact it is known \cite{j19} that among radial deformations ${\bf y}(\x) = r(|\x|)\frac{\x}{|\x|}$ the minimizer exists and satisfies $r(0)>0$, so that a cavity is formed at the origin. Thus enlarging the function space allows for the description of cavitation, a standard failure mechanism for rubber. Note that the growth rate of $\psi(\A)=|\A|^2+h(\det\A)$ for large $|\A|$ is what allows this description of cavitation. Had $|\A|^2$ been replaced by $|\A|^p$ for some $p>3$ then the model could not have predicted cavitation because every mapping in $W^{1,p}(\om;\R^3)$ is continuous by the Sobolev embedding theorem.

Suppose we use the prescription that the function space should be chosen to be the largest one for which the total free energy is defined.  At first sight this would seem to be the largest Sobolev space  $W^{1,1}(\om;\R^3)$.  However mappings $\y\in W^{1,1}(\om;\R^3)$ cannot be discontinuous across surfaces, so that to describe   cracks  an even larger space is needed. In the Francfort-Marigo formulation of the Griffith theory of brittle fracture \cite{bourdinetal08,francfortmarigo98} the space $SBV(\om;\R^3)$ is used instead, with the modified functional 
\be \label{elasticity1}I(\y)=\int_\om \psi(\nabla\y(\x))\,d\x+ \kappa{\mathcal H}^2(S_\y),\ee
where $\kappa>0$ and ${\mathcal H}^2(S_\y)$ is the two-dimensional Hausdorff measure  of the jump set $S_\y$, that gives the location of cracks in the reference configuration. More general surface energy terms can be considered to describe different kinds of fracture, such as 
\begin{equation}\label{elasticity2}\int_{S_\y}g(\y_+-\y_-,\nnu)\,d{\mathcal H}^2,\end{equation}
where $\nnu$ is the normal to $S_\y$ and $\y_+, \y_-$ are the approximate one-sided limits  of $\y$ on each side of the jump set.

Minimizers of \eqref{elasticity} may also be continuous but have jump discontinuities in $\nabla\y$ across surfaces. This occurs for models of elastic crystals, in which the surfaces across which $\nabla\y$ jumps represent phase boundaries (see \cite{j32,j40,bhattacharya03a}). These surfaces also carry energy, and models \cite{j58} have been proposed in the spirit of \eqref{elasticity1} in which $\nabla\y\in SBV(\om;M^{3\times 3})$ and sharp interfaces are allowed to compete energetically with diffuse ones.

As we will see, all these issues have their counterparts for liquid crystals, the important conclusion being that the function space is part of the model.

\section{Description of  defects}
 \label{defects}

We first consider the Oseen-Frank model for nematics with energy given by \eqref{oseenfrank} and $q_0=0$. Clearly $W$ given by \eqref{ofw} satisfies $W(\n,\nabla\n)\leq C|\nabla\n|^2$ for some $C>0$ whenever $|\n|=1$, and under the Ericksen inequalities \cite{ericksen66a,virgabook} 
\be \label{ericksenineq}
2K_1>K_2+K_4,\; \;K_2>|K_4|,\;\; K_3>0\ee
 we also have that $C'|\nabla\n|^2\leq W(\n,\nabla\n)$ for some constant $C'>0$. Hence a natural function space to use is $W^{1,2}(\om;S^2)$. Then the hedgehog $\tilde\n$ given by \eqref{hedgehog} is a solution of the Euler-Lagrange equations for $I_{OF}$ (interpreted in a suitable weak sense). In the one-constant case $K_1=K_2=K_3=K>0, K_4=0$ the functional reduces to the form
$$I(\n)=K\int_\om |\nabla\n|^2d\x,$$
for which the corresponding solutions to the Euler-Lagrange equation are known as {\it harmonic maps}. In this case it was proved by Brezis, Coron \& Lieb \cite{breziscoronlieb} that $\tilde\n$ is the absolute minimizer of $I$ subject to its own boundary conditions on $\partial\om$. Furthermore, by a result of Schoen \& Uhlenbeck \cite{schoenuhlenbeck,schoenuhlenbeck1} any minimizer $\n$ of $I$ with a given boundary condition is smooth except at a finite number of point defects located at points $\x(i)\in \om$, and it is shown in \cite{breziscoronlieb} that $$\n(\x)\sim \pm \Rz(i)\frac{\x-\x(i)}{|\x-\x(i)|}\mbox{ as }\x\to \x(i)$$ for some rotation $\R(i)\in SO(3)$. For the case of general elastic constants less is known. For $W$ satisfying \eqref{ericksenineq} Hardt, Lin \& Kinderlehrer \cite{hardtkinderlehrerlin} (for a different proof see Hong \cite{hong04}) showed that any minimizer is smooth except for a closed set $E$ of   Hausdorff dimension strictly less than 1, so that there are no line defects, but it is not known whether $E$ consists of just finitely many points. Another difference is that, as shown by Helein \cite{helein}, the hedgehog $\tilde\n$ is not a minimizer if $8(K_1-K_2)>K_3$, this condition being sharp \cite{cohentaylor}.

Now consider the line defect or disclination described by
\be 
\label{disclination}
\hat\n(\x)=(\frac{x_1}{r},\frac{x_2}{r},0)\;\;r=\sqrt{x_1^2+x_2^2}.
\ee
This is also a smooth solution of the Euler-Lagrange equation for $I_{OF}$ for $r\neq 0$. Since 
$|\nabla\hat\n(\x)|^2=\frac{1}{r^2}$ we have that $\hat\n\in W^{1,p}(\om;S^2)$ if and only if $1\leq p<2$, so that $I_{OF}(\hat\n)=\infty$. That the Oseen-Frank theory does not allow line defects with finite energy, whereas such defects are observed, has long been recognized as a serious deficiency. Various remedies have been proposed, two of them being the Landau - de Gennes and Ericksen theories. However, it is worth asking whether the deficiency can be remedied without introducing a higher-dimensional order parameter. There is an obvious way to do this, once we recognize that the infinite energy arises from the quadratic growth of $W(\n,\nabla \n)$ for {\it large} $|\nabla\n|$ while the quadratic character of $W$ suggests that it is a valid  approximation for {\it small} $|\nabla \n|$. So a possible solution to the deficiency is to change the growth rate of $W$ to be subquadratic as $|\nabla\n|\to \infty$ while retaining essentially the same form as \eqref{ofw}  for small $|\nabla\n|$. One of many ways to do this would be to replace $W(\n,\nabla\n)$ for fixed $p$, $1<p<2$, by
\be\label{neww}
W_\alpha(\n,\nabla\n)=\frac{2}{p\alpha}\left(\left(1+\alpha W(\n,\nabla \n)\right)^{\frac{p}{2}}-1\right),
\ee
where   $\alpha>0$ is small. Then $W_\alpha(\n,\nabla\n)\to W(\n,\nabla\n)$ as $\alpha\to 0$. Also, assuming the Ericksen inequalities \eqref{ericksenineq}, $W_\alpha$ satisfies the growth conditions 
\be 
\label{newgrowth}
C'_\alpha(|\nabla \n|^p-1)\leq W_\alpha(\n,\nabla\n)\leq C_\alpha|\nabla\n|^p,
\ee
for positive constants $C_\alpha, C'_\alpha$. Setting 
$$I_\alpha(\n)=\int_\om W_\alpha(\n,\nabla\n)\,d\x,$$
 we obtain from the right-hand inequality in \eqref{newgrowth} that $I_\alpha(\hat\n)<\infty$ as desired. Also $W_\alpha(\n,\D)$ is convex in $\D$. Indeed, fixing $\n$, the fact that  $W(\n,\cdot)$ is a positive quadratic form implies  that $\vertiii{\D} =W^\frac{1}{2}(\n,\D)$ defines a  norm on $M^{3\times 3}$. Therefore if $0\leq\lambda\leq 1$
\begin{eqnarray*}
&&\hspace{-.4in}W_\alpha(\n,\lambda\D+(1-\lambda)\D')\\
&=&\frac{2}{p\alpha}\left(\left(\left(1+\alpha\vertiii{\lambda\D+(1-\lambda)\D'}^2\right)^\frac{1}{2}\right)^p-1\right)\\
&\leq &\frac{2}{p\alpha}\left(\left(\lambda\left(1+\alpha\vertiii{\D}^2\right)^\frac{1}{2}+\left(1-\lambda\right)\left(1+\alpha\vertiii{\D'}^2\right)^\frac{1}{2}\right)^p-1\right)\\
&\leq &\frac{2}{p\alpha}\left( \lambda\left(\left(1+\alpha\vertiii{\D}^2\right)^\frac{p}{2}-1\right)+\left(1-\lambda\right)\left(\left(1+\alpha\vertiii{\D'}^2\right)^\frac{p}{2}-1\right)\right)\\
&=&\lambda W_\alpha(\n,\D)+(1-\lambda)W_\alpha(\n,\D'),
\end{eqnarray*}
where we have used the convexity of $(1+\alpha\vertiii{\cdot}^2)^\frac{1}{2}$ and of $|\cdot|^p$. This convexity, together with the left-hand inequality in \eqref{newgrowth}, is enough to ensure, by standard methods of the calculus of variations, that  there exists at least one minimizer of $I_\alpha(\n)$ satisfying given suitable boundary conditions. 
Of course the particular replacement \eqref{neww} doesn't have a physical motivation, but it indicates the possibilities. 

The use of a modified free-energy density with subquadratic growth for large $|\nabla\n|$, such as $W_\alpha(\n,\nabla\n)$, can also widen the class of possible boundary conditions. If $\om\subset \R^3$ has smooth boundary and a sufficiently smooth 
unit vector field $\bf N$ is given on the boundary $\partial \om$, then
it is known \cite{hardtlin} that there is a unit vector field
$\n\in W^{1,2}(\om;S^2)$ with $\n={\bf N}$ on $\partial\om$. However, if, for example, $\om=(0,1)^3$ is a cube and $\bf N$ is the 
inward normal to the boundary, then \cite{bedfordthesis} there is no such $\n$.
Thus the Oseen-Frank theory does not apply to homeotropic boundary conditions on a cube,
although a  theory with subquadratic growth would allow  such boundary conditions.

\begin{figure}[ht]
\label{orien}
\centering
\includegraphics[scale =0.37]{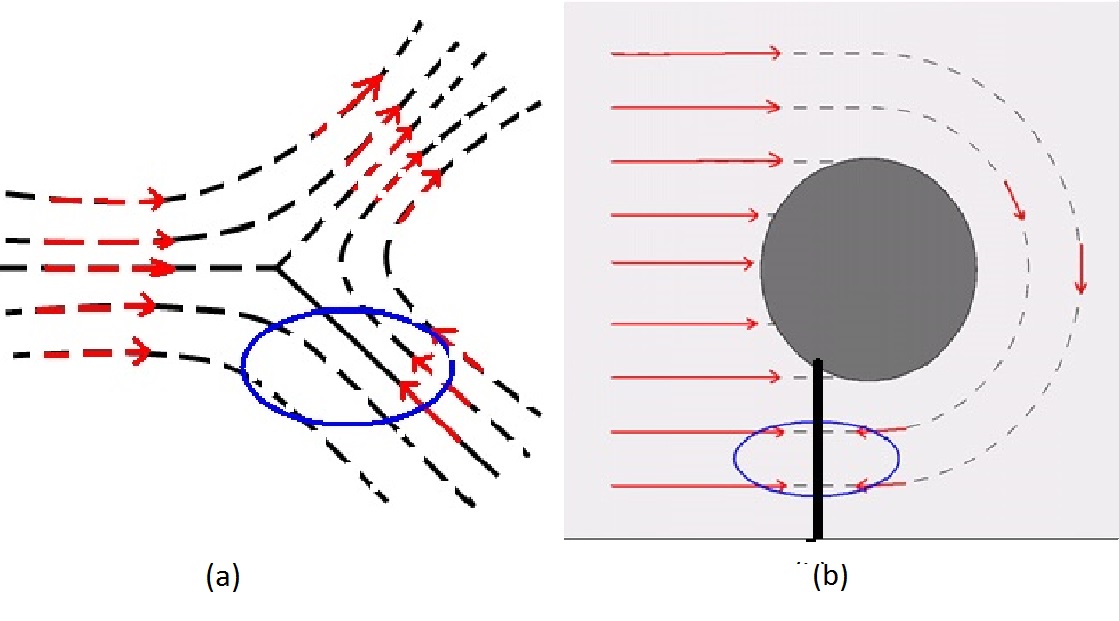}
\caption{Non-orientable line fields: (a) index $\frac{1}{2}$ singularity,   (b) smooth line field in non simply-connected region outside a cylinder. 
In both diagrams the line field is tangent to the dashed lines, with the perpendicular component zero. Attempting to orient the line field in a Sobolev space gives a contradiction in the region inside the elliptical cylinder. By introducing a discontinuity along the surfaces indicated by the solid lines, the line fields can be represented by a vector field in $SBV$.}
\end{figure}

However, the above idea is not in itself sufficient for line singularities to be describable within a pure director theory because of the existence of observed  non-orientable line defects, such as the index $\frac{1}{2}$ defect illustrated in Fig.~\ref{orien} (a). The issue of when line fields are orientable was discussed in \cite{j61}, where it was shown that if $\om$ is simply-connected   any line field $\n\otimes\n\in W^{1,2}(\om;M^{3\times 3})$ is {\it orientable}; that is there is   a unit vector field $\m\in W^{1,2}(\om;S^2)$ with $\m(\x)=\kappa(\x) \n(\x)$ with $\kappa(\x)=\pm 1$ (equivalently $\m(\x)\otimes \m(\x)=\n(\x)\otimes\n(\x)$) for almost every $\x\in\om$. In non simply-connected domains even smooth line fields need not be orientable (see Fig. \ref{orien} (b)). For the index $\frac{1}{2}$ singularity the result in \cite{j61} does not apply, for although $\om$ can be assumed simply-connected  the corresponding line field does not belong to $W^{1,2}(\om;M^{3\times 3})$. 
If, however, we allow $\n$ to reverse orientation across suitable surfaces, as illustrated in Fig. \ref{orien}, we can represent such non-orientable line fields by vector fields that belong to $SBV(\om;S^2)$. A result in this direction is the following:
\begin{thm}[Bedford \cite{bedfordcholesterics}]
\label{SBVorient}
  Let $\Q=s\left(\n\otimes\n-\frac{1}{3}
{\bf 1}\right)\in W^{1,2}(\om;M^{3\times 3})$, where
$s\neq 0$ is constant. Then there exists a unit vector field ${\bf m}\in SBV(\om;S^2)$
such that $\m\otimes \m=\n\otimes \n$, and if $\x\in S_\m$ then $\m_+(\x)=-\m_-(\x)$.
\end{thm} 
This theorem applies to the case illustrated in Fig. \ref{orien} (b), but not to the index $\half$ singularity in Fig. \ref{orien} (a), for which an extension to $\Q\in W^{1,p}(\om;M^{3\times 3})$ for $1<p<2$ would be required.
It provides a first possible use of director fields $\n$ that jump across surfaces, though it is as  a mathematical device rather than representing an actual physical discontinuity.

For the Ericksen theory  we can model point  defects and 
disclinations by finite energy  configurations in which $\n$ is
 discontinuous and  $s=0$ at the defect (melting core).  In this case there is no need to change the growth rate of the free-energy density at infinity. For example,
if we consider the special case when 
$$I_{LdG}(\Q)=\int_\om[K|\nabla\Q|^2+\psi_B(\Q)]\,d\x,$$
then the ansatz
\be 
\label{uniaxials}
\Q(\x)=s(\x)\left(\n(\x)\otimes\n(\x)-\frac{1}{3}{\bf 1}\right)
\ee
gives the functional
$$I_E(s,\n)=\int_\om [K(|\nabla s|^2+2s^2|\nabla\n|^2) +\psi_B(s)]\,ds,$$
where $\psi_b(s)=\hat\psi(\frac{2s^2}{3}, \frac{2s^3}{27})$ (see \eqref{isotropy}).
Then $\n$ can have a singularity at a point or curve which has finite energy 
because $s$ can tend to zero sufficiently fast as the point or curve is
approached to make $I_E(s,\n)$ finite. However for non simply-connected domains or index $\frac{1}{2}$ defects there is the same
 orientability
problem as in the Oseen-Frank theory, which can be `cured' by allowing jumps from $\n$ to $-\n$
 across surfaces. Here we need to remove the closure of the set of points where $s=0$ in order to assert the existence of a corresponding $SBV$ mapping.
\begin{thm} [Bedford \cite{bedfordcholesterics}]  Let $\Q\in W^{1,2}(\om;M^{3\times 3})$ have the uniaxial form \eqref{uniaxials} with $s\in C(\Omega)$. 
Then $s\in W^{1,2}(\om)$ and there
exists a  vector field 
 ${\bf m}\in SBV_{\rm loc}(\om\setminus {\{s=0\}};S^2)$ such 
that ${\bf m}\otimes{\bf m}=\n\otimes \n$ .
\end{thm}
There is also the possibility of `genuine' planar defects in the Ericksen
theory (see Section \ref{planar}).

In the Landau - de Gennes theory minimizers $\Q$ of the free-energy functional $I_{LdG}(\Q)$
are expected to be smooth \cite{gartlanddavis} because the 
Euler-Lagrange equation is elliptic. (There is, however, an issue if the bulk energy satisfies the blow-up condition \eqref{blowup}, when it is not obvious in general that the minimum eigenvalue $\lambda_{\rm min}(\Q(\x))$ is bounded away from zero, which makes it unclear whether minimizers satisfy the Euler-Lagrange equation. This property of the eigenvalues is proved in the one-constant case in \cite{u9,j59}, but otherwise is an open problem.)
Therefore in the Landau - de Gennes theory defects are not described by singularities of $\Q$, but rather by 
regions where $|\nabla\Q|$ is very high, in which a continuous choice 
of eigenvectors for $\Q$ may not possible.

The situation might be different for free-energy densities $\psi(\Q,\nabla\Q)$
which are convex (or quasiconvex) but not quadratic in $\nabla \Q$. For such integrands there is a counterexample of
\v{S}ver\'{a}k \& Yan \cite{sverakyan00} which has a singular minimizer of the form
$$\Q(\x)=|\x|\left(\frac{\x}{|\x|}\otimes\frac{\x}{|\x|}-\frac{1}{3}{\bf 1}\right).$$

\section{Planar defects}
\label{planar}
In this section we explore situations in which it may be physically reasonable to allow order parameters to jump across surfaces. In this context it is natural to  consider an energy functional analogous to \eqref{elasticity1},\eqref{elasticity2} given by 
\begin{equation}\label{sbvlc}
I(\n)=\int_\om W(\n,\nabla\n)\,d\x + \int_{S_\n}f(\n_+,\n_-,\nnu)\,d{\mathcal H}^2,
\end{equation}
for $\n\in SBV(\om;S^2)$, where $\nnu$ is the normal to the jump set $S_\n$. Here $W(\n,\nabla\n)$ is assumed to have the Oseen-Frank form \eqref{ofw} or be modified so as to have subquadratic growth as suggested in Section \ref{defects}. Such  models do not seem to have been previously considered in the liquid crystal literature, although  a similar model is mentioned in \cite[Section 4.6.4]{ambrosioetal00}. 

We will suppose that $f:S^2\times S^2\times S^2\to[0,\infty)$ is a continuous function which is  frame-indifferent, i.e.
\be 
\label{frind}
f(\Rz\n_+,\Rz\n_-,\Rz\nnu)=f(\n_+,\n_-,\nnu) \mbox{ for all } \Rz\in SO(3), \n_+,\n_-,\nnu\in S^2,
\ee 
and that $f$ is invariant to reversing the signs of $\n_+, \n_-$, reflecting the statistical head-to-tail symmetry of nematic and cholesteric molecules, so that
\be 
\label{headtotail}
f(-\n_+,\n_-, \nnu)=f(\n_+,-\n_-,\nnu)=f(\n_+,\n_-,\nnu).
\ee
(Note that taking $\Rz={\rm diag}\,(1,-1,-1)$ in \eqref{frind} we have that $f(\n_+,-\n_-,-\nnu)=f(\n_+,\n_-,\nnu)$, so that by \eqref{headtotail} we also have that $f(\n_+,\n_-,-\nnu)=f(\n_+,\n_-,\nnu)$.)

 It is proved in 
 \cite{u14} 
that a necessary and sufficient condition for  a continuous function $f:S^2\times S^2\times S^2\to[0,\infty)$ to satisfy   \eqref{frind} and \eqref{headtotail} is that 
\be 
\label{represent}
f(\n_+,\n_-,\nnu)=g((\n_+\cdot\n_-)^2, (\n_+\cdot\nnu)^2,(\n_-\cdot\nnu)^2,(\n_+\cdot\n_-)(\n_+\cdot\nnu)(\n_-\cdot\nnu))
\ee 
for a continuous function $g:D\to[0,\infty)$, where
$$D=\{(\alpha,\beta,\gamma,\delta):\alpha,\beta,\gamma\in [0,1], \delta^2=\alpha\beta\gamma, \alpha+\beta+\gamma-2\delta\leq 1\}.$$
Since $D$ is compact, the function $g$ in the theorem can be extended to a continuous function on $\R^4$.
 
It is also natural to assume that 
\be 
\label{zeroenergy}f(\az,\az,\nnu)=0\mbox{ for all }\az,\nnu,
\ee
 so that the energy of a jump tends to zero as $|\n_+-\n_-|\to 0$. In the case when $f$ does not depend on $\nnu$ we thus have that $f(\n_+,\n_-)=g((\n_+\cdot\n_-)^2)$ with $g(1)=0$, a simple  example being 
\be 
\label{fex}
f(\n_+,\n_-)=k'(1-(\n_+\cdot\n_-)^2)^\frac{r}{2}
\ee
for  constants $k'>0,\; 0<r< 1$. 

In fact it is not clear how to analyze the functional \eqref{sbvlc} directly, due to the difficulty of controlling changes in the sign of $\n$.  Instead we consider in Sections \ref{order} and \ref{smectics} below the corresponding functionals 
\be 
\label{sbvqlc}
I(\Q)=\int_\om \psi(\Q,\nabla\Q)\,d\x +\int_{S_\Q}F(\Q_+,\Q_-,\nu)\,d{\mathcal H}^2,
\ee
where $\Q\in SBV(\om;M^{3\times 3})$ satisfies the uniaxial constraint 
\be 
\label{uniaxiala}
\Q=s\left(\n\otimes\n-\frac{1}{3}{\bf 1}\right)
\ee
 for fixed $s>0$.  Then \eqref{represent} takes the form
\be 
\label{representQ}
F(\Q_+,\Q_-,\nu)=G(\Q_+\cdot\Q_-, \Q_+\nu\cdot\nu, \Q_-\nu\cdot\nu, \Q_+\nu\cdot\Q_-\nu)
\ee
for some function $G$, corresponding to the condition that $F$ is an isotropic function of $\Q_+,\Q_-$ and $\nu\otimes\nu$ (see
\cite[p910, first line of (2.37)]{smith71}).  The special case \eqref{fex} takes the form
\be 
\label{fexQ}
F(\Q_+,\Q_-)=k|\Q_+-\Q_-|^r,
\ee 
where $k=2^{-\frac{r}{2}}s^{-r}k'$.
The condition   $0<r<1$ is needed to ensure existence of a minimizer for \eqref{sbvlc}   (see Section \ref{smectics} below). 

We remark that models of the form \eqref{sbvlc}, \eqref{sbvqlc} can potentially give a description of line defects of finite energy even when $W, \psi$ are quadratic in $\nabla\n, \nabla\Q$ respectively, since, for example, we can  surround a line defect with a cylindrical core inside which we choose $\n$ or $\Q$ to be constant. This gives   a configuration  close to that of a line defect having finite total energy, namely the elastic energy outside the core plus the interfacial energy of the core boundary.

\subsection{Nematic elastomers}
\label{nematicelastomers}
In Section \ref{solids} we discussed models of martensitic phase transformations, in which the deformation gradient can jump across surfaces representing phase boundaries, which can be atomistically sharp or diffuse (but still of thickness only a few atomic spacings). Typical martensitic microstructures  exhibit fine laminates formed by parallel such phase boundaries separating two different variants of the low temperature martensitic phase. Nematic elastomers (see \cite{warnerterentjev}) exhibit similar laminated structures \cite{kundlerfinkelmann}. The energy functional for nematic elastomers proposed in \cite{bladonetal} is given by
\be 
\label{elastomer}
I(\y,\n)=\int_\om \frac{\mu}{2}\left(D\y (D\y)^T\cdot L_{a, \n}^{-1}-3\right)\,d\x,
\ee
where
\be 
\label{Ltensor}
L_{a, \n}=a^\frac{2}{3}\n\otimes\n+a^{-\frac{1}{6}}({\bf 1}-\n\otimes\n)
\ee
and $\mu>0, a>0$ are material parameters. The material is assumed incompressible, so that $\y$ is subjected to the constraint $\det \nabla \y=1$.
By minimizing 
 the integrand over $\n\in S^2$ we obtain the purely elastic energy
\be 
\label{elasticenergy}
I(\y)=\int_\om W(\nabla\y)\,d\x,
\ee
where
\be 
\label{rel}
W(\A)=\frac{\mu}{2}\left(a^{-\frac{2}{3}}v_1^2(\A)+a^\frac{1}{3}(v_2^2(\A)+v_3^2(\A))\right),
\ee
and $v_1(\A)\geq v_2(\A)\geq v_3(\A)>0$ denote the singular values of $\A$, that is the eigenvalues of $\sqrt{\A^T\A}$.

As discussed by De Simone \& Dolzmann \cite{desimonedolzmann}, the free-energy function  \eqref{rel} is not quasiconvex, and admits minimizers in which $\nabla \y$ jumps across planar interfaces, so that the minimizing $\n$ of the integrand also jumps. Of course the functional \eqref{elastomer} ignores Frank elasticity, i.e. terms in $\nabla\n$, but the experimental observations might suggest that even with such terms allowing jumps in $\n$ may be a useful approximation.

\subsection{Order reconstruction problems}
\label{order}
 A number of authors (for example \cite{barberobarberi, carboneetal, palffyetal, zannonithinfilm}) have considered situations in which a nematic liquid crystal is confined between two plates or surfaces which are very close together, and on which the orientations of the director are different. Similar situations arise in thin films of nematic liquid crystals, in which the free surface replaces one of the plates. Below a critical separation distance between the plates the system prefers a biaxial order in which the director may change discontinuously due to an exchange of eigenvalues. In the theory of Palffy-Muhoray, Gartland \& Kelly \cite{palffyetal} the Landau - de Gennes theory is used with energy functional
 \be 
\label{ldg1}
I_{LdG}(\Q)=\int_{\om_\delta} [\psi_B(\Q)+\psi_E(\Q,\nabla\Q)] \,d\x,
\ee
with $\psi_B$ having the quartic form \eqref{bulk} and $\psi_E=K|\nabla\Q|^2$, with the liquid crystal occupying the region $\Omega_\delta=(0,1)^2\times (0,\delta)$  lying between two plates at $x_3=0$ and $x_3=\delta>0$. We will consider the energy \eqref{ldg1} with  $\psi_B=\psi_B(\Q)$ everywhere defined, continuous and bounded below (such as  the quartic form \eqref{bulk}) or defined and continuous for $\lambda_{\rm min}(\Q)>-\frac{1}{3}$ and satisfying \eqref{blowup}, and with $\psi_E=\psi_E(\nabla\Q)$ a positive quadratic form in $\nabla\Q$ (in particular $L_4=L_5=0$ in \eqref{elastic}).  We will assume that the boundary conditions satisfied on the plates $x_3=0$ and $x_3=\delta$ are planar and homeotropic respectively, so that
\be 
\label{qbc}
\Q(x_1,x_2,0)=\Q^{(0)},  \;\;\Q(x_1,x_2,\delta)=\Q^{(1)},\mbox{ for all }(x_1,x_2)\in (0,1)^2, \\
\mbox{where }\Q^{(0)}=s_1\left(\e_1\otimes\e_1-\frac{1}{3}{\bf 1}\right),\;\; \Q^{(1)}:=s_2\left(\e_3\otimes\e_3-\frac{1}{3}{\bf 1}\right),\nonumber
\ee
and $0<s_1\leq 1$, $0<s_2\leq 1$ (in \cite{palffyetal} it is assumed that $s_1=s_2$). For lateral boundary conditions we assume for simplicity that $\Q$ is periodic in $x_1,x_2$, so that
\be 
\label{latbc}
\Q(0,x_2,x_3)=\Q(1,x_2,x_3),\;\;\Q(x_1,0,x_3)=\Q(x_1,1,x_3),
\ee
 for all $(x_2,x_3)\in(0,1)\times(0,\delta)$ and $(x_1,x_3)\in (0,1)\times (0,\delta)$ respectively. We want to study the behaviour of minimizers of \eqref{ldg1} subject to \eqref{qbc}, \eqref{latbc} in the limit when $\delta\to 0$. To this end it is convenient to rescale $\Omega$ by a factor $\frac{1}{\delta}$ in the $x_3$-direction and set
\be 
\label{rescale}
\Pz(x_1,x_2,x_3)=\Q(x_1,x_2,\delta x_3),
\ee
so that $I_{LdG}(\Q)=\delta^{-1}E^\delta(\Pz)$, where
\be 
\label{rescaleI}
E^\delta(\Pz)=\int_C [\delta^2\psi_B(\Pz)+\psi_E(\delta\Pz_{,1},\delta\Pz_{,2},\Pz_{,3})]\,d\x
\ee 
and $C=(0,1)^3$ is the unit cube. The following theorem shows that in the limit $\delta\to 0$ minimizers become linear in $x_3$ (consistent with \cite[Fig.1]{palffyetal}, where $\delta$ is small).
\begin{thm}
\label{linear}
Let $\Pz^\delta$ be a minimizer of $E^\delta$ subject to the rescaled boundary conditions $\Pz^\delta(x_1,x_2,0)=\Q^{(0)}, \Pz^\delta(x_1,x_2,1)=\Q^{(1)}$ and $\Pz^\delta$ periodic in $x_1,x_2$. Then as $\delta\to 0$
$$\Pz^\delta\to \bar\Pz,\; \Pz^\delta_{,3}\to \bar\Pz_{,3}, \;\delta \Pz^\delta_{,1}\to 0,\; \delta \Pz^\delta_{,2}\to 0\; \mbox{ in }L^2(C),$$
 $($i.e. for example, $\int_C|\Pz^\delta-\bar\Pz|^2\,dx\to 0),$
  where
$$\bar\Pz(x_1,x_2,x_3)= (1-x_3)\Q^{(0)}+x_3\Q^{(1)}.$$
\end{thm}
Theorem \ref{linear} shows that, for sufficiently small $\delta$, $\Q$ is given approximately by
$$ \Q(\x)= (1-\delta^{-1}x_3)\Q^{(0)}+ \delta^{-1}x_3\Q^{(1)},$$
for which the maximum eigenvalue is given by 
\be 
\label{s}
\lambda_{\rm max}(\x)=\left\{\begin{array}{ll}\frac{2}{3}s_1-\frac{2s_1+s_2}{3}\frac{x_3}{\delta}&\mbox{if } 0\leq x_3\leq\frac{s_1}{s_1+s_2}\delta \\
-\frac{1}{3}s_1+\frac{2s_2+s_1}{3}\frac{x_3}{\delta}&\mbox{if }  \frac{s_1}{s_1+s_2}\delta\leq x_3\leq 1,
\end{array}\right.
\ee
with corresponding eigenvector
\be 
\label{evector}
\n(\x)=\left\{\begin{array}{ll}\e_1&\mbox{if } 0\leq x_3\leq\frac{s_1}{s_1+s_2}\delta \\ \e_3&\mbox{if } \frac{s_1}{s_1+s_2}\delta\leq x_3\leq 1.
\end{array}\right.
\ee
Thus the director $\n$ has a discontinuity on the plane $x_3=\frac{s_1}{s_1+s_2}\delta$. In \cite{palffyetal} it is shown numerically that a bend solution corresponding for large $\delta$ to that of the Oseen-Frank theory (given by $\n(\x)=\cos\left(\frac{\pi x_3}{2\delta}\right)\e_1+\sin\left(\frac{\pi x_3}{2\delta}\right)\e_3$) bifurcates from the small $\delta$ minimizer at a critical value of $\delta$. Similar numerical simulations are presented in \cite{bisietal03}, while Lamy \cite{lamy14} gives rigorous results establishing such a bifurcation in a relevant parameter regime.

Of course it is not immediately clear if the Landau - de Gennes theory applies at the very small length scales (of the order of 10nm) at which it predicts an eigenvalue exchange. However, in \cite{carboneetal} a good correlation between experimentally measured forces and simulations based on the Landau - de Gennes model is found, and a number of similar studies referred to which draw a similar conclusion. (The simulations reported in \cite{carboneetal} are based on the methodology in \cite{lombardoetal}, which uses a dynamical minimization based on the quartic bulk energy \eqref{bulk} and an elastic energy of the form \eqref{elastic} with three nonzero coefficients $L_1, L_3, L_4$; without a constraint on the eigenvalues of $\Q$, such as provided by \eqref{blowup}, the corresponding total free energy \eqref{ldg} is unbounded below for any boundary conditions \cite{j59}, but this may not be a problem for appropriate dynamic simulations --- see \cite{iyeretal} for a discussion. In \cite{barberietal04}, on the other hand, there are just two nonzero elastic coefficients $L_1, L_3$ and the total energy is bounded below.)   

\begin{figure}[h] 
  \centering
  \includegraphics[width=5in,height=3.14in,keepaspectratio]{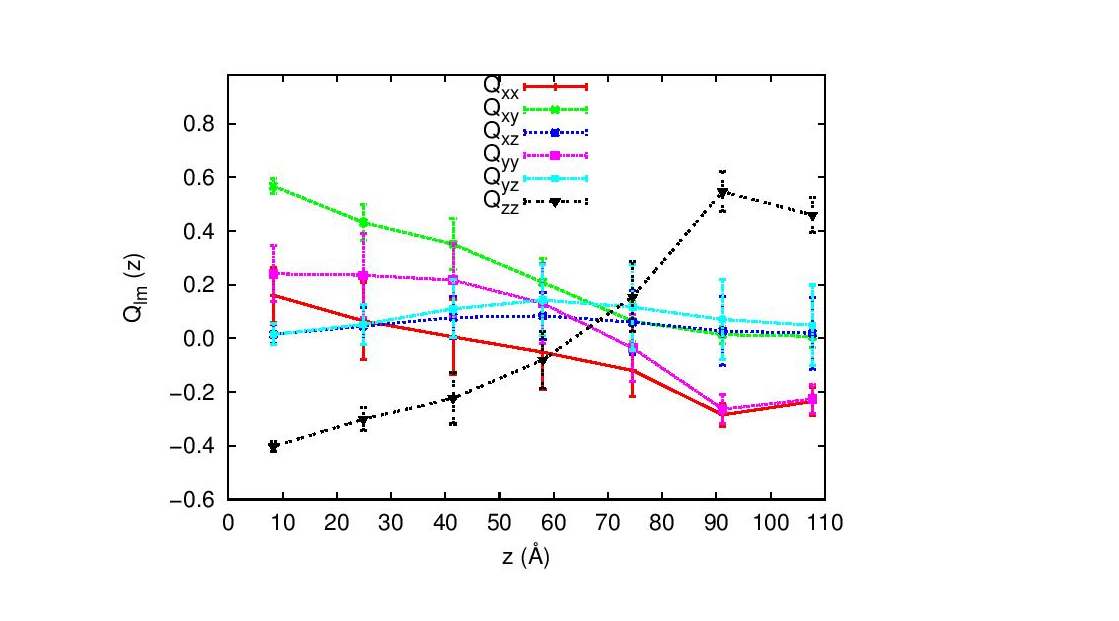}
  \caption{$\Q$-tensor values for 12nm thin film between silicon and a vacuum obtained by atomistic molecular dynamics \cite{zannonithinfilm}  (courtesy C. Zannoni and M. Ricci); the error bars correspond to the standard deviations obtained when averaging over the molecules in each 1nm thick slice and over the trajectory. }
  \label{Qfilm}
\end{figure}

\begin{figure}[h] 
  \centering
  \includegraphics[width=5.67in,height=1.9in,keepaspectratio]{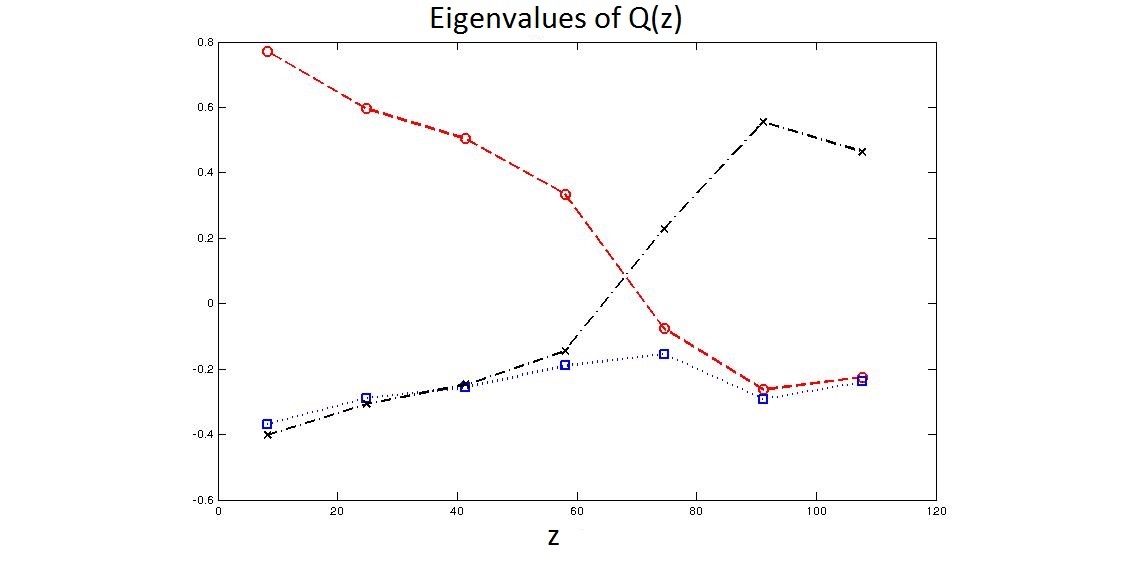}
  \caption{Eigenvalues of $\Q(z)$ as given in Fig.~\ref{Qfilm}}
  \label{qfilm1}
\end{figure}

In this context the recent atomistic molecular dynamics simulations \cite{zannonithinfilm} of a thin film of 5CB nematic liquid crystal on a hydrogenated silicon substrate are of interest. In the simulations the upper surface of the film was assumed to be in contact with a vacuum. The two cases of 12 and 24nm thick    films were studied. In Fig. \ref{Qfilm} 
the values of the $\Q$-tensor as a function of $z=x_3$ for the thinner film are shown (this figure is not included in \cite{zannonithinfilm}), while in Fig. \ref{qfilm1}
the corresponding values of the three eigenvalues of the mean value of $\Q$ are given\footnote{The graph of the largest eigenvalue has the same qualitative form as that in \cite[Fig. 4]{zannonithinfilm} but the values are less. This is because the values in \cite[Fig. 4]{zannonithinfilm} were obtained as the average   of the largest eigenvalue of $\Q(z)(jT)$ over the trajectory $j=1,\cdots,N$, where $\Q(z)(jT)$ denotes the $\Q$ matrix calculated at time $jT$ by averaging over the   molecules in the 1nm slice centred at $z$; since the largest eigenvalue of a symmetric  matrix is a convex function of the matrix the maximum eigenvalue of the average of $\Q(z)(jT)$ over $j$ is by Jensen's inequality less than or equal to the average over $j$  of the largest eigenvalue of $\Q(z)(jT)$.}. There is clear evidence of an eigenvalue exchange at about $z=6.81$nm, with the director (the eigenvector corresponding to the maximum eigenvalue) jumping from about $(0.5037,0.8466,-0.1721)$ to $(0.0108,0.1991,0.9799)$. 

In summary, there is good evidence based on experiment, theory and simulation of discontinuities in the director, that is in the eigenvector of $\Q$ corresponding to its maximum  eigenvalue, across surfaces. Whereas this can be understood using the Landau - de Gennes theory, the existence of a critical separation between the plates below which a discontinuity appears is also predicted by a director model of the type \eqref{sbvqlc}. To see why this is so, consider for simplicity the case when $s_1=s_2=s$ with the one constant approximation $\psi(\Q,\nabla\Q)=K|\nabla\Q|^2$, with $K>0$, and the simple form \eqref{fexQ} for the jump set energy, with the liquid crystal occupying as before the region $\om_\delta=(0,1)^2\times (0,\delta)$. Then the problem becomes that of minimizing
\be 
\label{disc}
I(\Q)=\int_{\om_\delta} K|\nabla\Q|^2d\x +k \int_{S_\Q} |\Q_+-\Q_-|^rd{\mathcal H}^2
\ee 
for $\Q\in SBV(\om:M^{3\times 3})$ satisfying the uniaxial constraint \eqref{uniaxiala} and the boundary conditions 
\be   
\label{qbc1}
\Q(x_1,x_2, 0)=s\left(\e_1\otimes\e_1-\frac{1}{3}\1\right),\;\; \Q(x_1,x_2,\delta)=s\left(\e_3\otimes\e_3-\frac{1}{3}\1\right) 
\ee 
for all $(x_1,x_2)\in (0,1)^2$, with $\Q$ periodic in $x_1,x_2$ (see \eqref{latbc}).  Some care is needed when interpreting the boundary conditions \eqref{qbc1} and periodicity, since it is possible that $\Q$ might jump at the boundary $\partial\Omega_\delta$ of $\om_\delta$. This is handled by minimizing $I(\Q)$ among $\Q\in SBV_{\rm loc}(\R^2\times (-1,\delta+1); M^{3\times 3})$ which satisfy $\Q(x_1,x_2,x_3)= s\left(\e_1\otimes\e_1-\frac{1}{3}\1\right)$ for $-1<x_3<0$, $\Q(x_1,x_2,x_3)=s\left(\e_3\otimes\e_3-\frac{1}{3}\1\right)$ for $\delta<x_3<\delta+1$, and $\Q(x_1+1,x_2,x_3)=\Q(x_1,x_2+1,x_3)=\Q(x_1,x_2,x_3)$ for all $(x_1,x_2,x_3)\in \R^2\times (-1,\delta+1)$. With this interpretation $S_\Q$ can be partly on  $\partial\Omega_\delta$.

Candidates for minimizers of $I$ are the two smooth $\Q$ given by
\be 
\label{Qsmooth}
\Q^\pm(\x)=\frac{s}{2}\left(\begin{array}{ccc}\frac{1}{3}+\cos\frac{\pi x_3}{\delta} & 0 & \pm \sin \frac{\pi x_3}{\delta}\\0 & -\frac{2}{3} & 0\\
\pm  \sin \frac{\pi x_3}{\delta}& 0 & \frac{1}{3} -\cos \frac{\pi x_3}{\delta}\end{array}\right),
\ee
which are the minimizers of $\int_{\om_\delta}|\nabla\Q|^2d\x$ among uniaxial $\Q\in W^{1,2}(\om_\delta;M^{3\times 3})$ satisfying the boundary conditions, and which correspond to the two Oseen-Frank solutions in which the line field rotates anticlockwise (resp. clockwise) in the $(x_1,x_3)$ plane from horizontal to vertical. Then we have the following result.
\begin{thm}[\cite{u14}]
\label{critical}For any $\delta>0$ there exists at least one minimizer $\Q\in SBV(\om_\delta:M^{3\times 3})$ of $I$ subject to the boundary conditions \eqref{qbc1} and periodicity in $x_1, x_2$. 
\end{thm}
We conjecture that   there is a small $\delta_0>0$ such that if $\delta>\delta_0$ then $\Q^\pm$ are the only minimizers, while if $0<\delta<\delta_0$ then   any minimizer $\Q$ has a single jump with jump set $S_\Q=\{\x:x_3=\gamma(\delta)\}$, where $0<\gamma(\delta)<\delta$. For results in this direction see \cite{u14}.
That $\Q^\pm$ are not minimisers for $\delta$ sufficiently small is easily seen. In fact, $|\nabla \Q^\pm|=\frac{C}{\delta}$ for some $C>0$, so that $I(\Q^\pm)= Ks^2\frac{C^2}{\delta}$. But if  
\be 
\label{Qhat}
\hat\Q(\x)=\left\{\begin{array}{ll}s\left(\e_1\otimes\e_1-\frac{1}{3}\1\right)&\mbox{if } 0<x_3<\frac{\delta}{2}\\ s\left(\e_3\otimes\e_3-\frac{1}{3}\1\right)&\mbox{if } \frac{\delta}{2}<x_3<\delta\end{array}\right.
\ee
then $I(\hat\Q)=ks^r2^\frac{r}{2}$, so that $I(\hat\Q)<I(\Q^\pm)$ if $s^{r-2}\frac{\delta k}{K}<2^{-\frac{r}{2}}C^2$. The numerical computations illustrated in \cite[Fig. 3]{bisietal03} do suggest quite a sharp transition across the eigenvalue exchange interface between predominantly uniaxial regions, indicating that a discontinuous director model is not so unreasonable.

\subsection{Smectics}
\label{smectics}
\begin{figure}[h] 
  \centering
  \includegraphics[width=4.85in,height=1.21in,keepaspectratio]{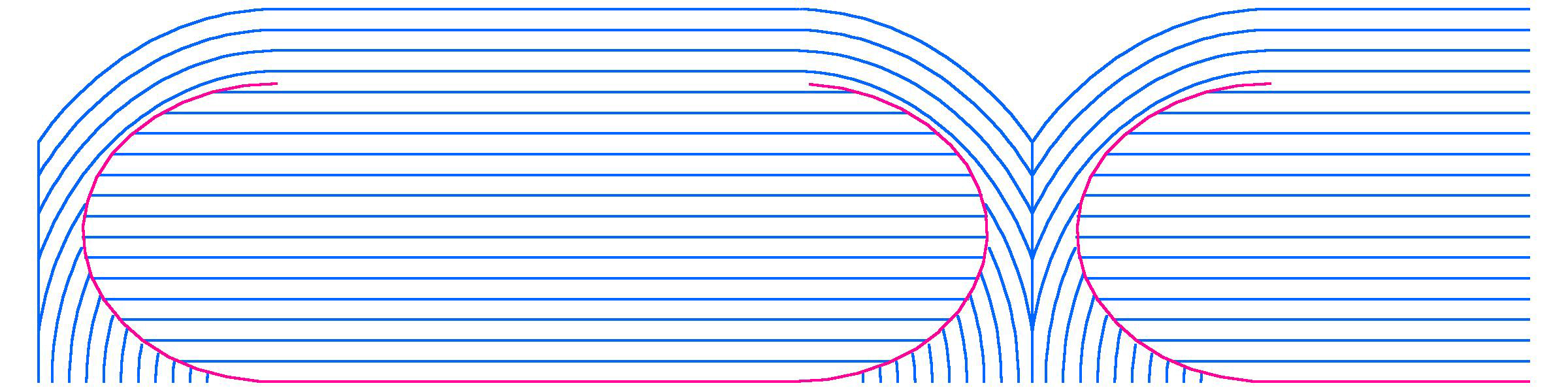}
  \caption{Smectic layer structure in 8CB thin film as identified in \cite{michel06} (courtesy Emmanuelle Lacaze).}
  \label{flat}
\end{figure}

Striking examples of apparent surface discontinuities in the director arise in recent experiments of E. Lacaze and coauthors  on  smectic-A liquid crystals \cite{lacazemicheletal,michel04, michel06,zapponelacaze,zapponeetal}. In these experiments a thin film of 8CB smectic A liquid crystal was deposited on a crystalline substrate with the upper  surface of the film in contact with air. The film thicknesses varied from 50nm to 1-2$\mu$m, so that the thinnest film was considerably thicker than those   in the nematic simulations of \cite{zannonithinfilm}. In \cite{lacazemicheletal,michel04,michel06} the substrate was ${\rm MoS}_2$, while in \cite{zapponeetal,zapponelacaze} it was crystalline mica. In \cite{michel04} X-ray diffraction measurements are reported that are consistent with the 8CB molecules at the substrate lying parallel to certain directions in the plane of the substrate, with the smectic layers at the substrate being perpendicular to the substrate. At the interface between the liquid crystal and the air, however, the molecules lie perpendicular to the interface, and the smectic layers are parallel to it. Thus, as in the order reconstruction problem, there are antagonistic boundary conditions at the two interfaces, leading to interesting defect structures. The structures identified in \cite{michel06} on the basis of X-ray diffraction measurements consist of a periodic family of flattened hemicylinders separated by vertical walls, as shown in cross-section in Fig. \ref{flat}.  
The curved and vertical interfaces are of the type called `curvature walls', for which various theoretical descriptions have been proposed  \cite{blanckleman99,degennesprost,dozovdurand94,williamskleman}. For thicker films these walls are not generally observed, and instead arrays of focal-conic defects  with vertical disclination lines are typical (see, for example, \cite{zapponeetal12} and for the case of degenerate planar anchoring \cite{designolle2006afm}). Two features of the defect morphology in Fig.~\ref{flat} are worth noting. First, equally spaced smectic layers can only match perfectly at a  curvature wall  if the tangent plane to the wall bisects the smectic planes on either side; otherwise there will be a mismatch leading to `dislocations' at the wall. Second, it is not obvious how the planar boundary conditions on the substrate are satisfied in the thin region at the bottoms of the squashed cylinders. However, in a later experiment \cite{coursaultetal2014} (see especially Fig.~2) with a different noncrystalline PVA substrate, a morphology for this region satisfying planar boundary conditions was identified. A somewhat different structure for the mica substrate was proposed in \cite[Fig.~7]{zapponeetal}. 

 Following de Gennes \cite{degennes72} it has become usual to model smectics by means of a complex order parameter $\Psi(\x)=r(\x) e^{i\phi(\x)}$, in terms of which the molecular density is given by 
$$\rho_0+\rho(\x)=\rho_0+{\rm Re}\Psi(\x)=\rho_0+r(\x)\cos\phi(\x),$$
where $\rho_0>0$ is a constant average density. Thus $\rho(\x)$ describes the fluctuations in the density due to the smectic layers, and $\nabla\phi$ gives the normals to the layers.   Various free-energy densities for smectics have been proposed (see, for example, \cite{chenlubensky,hanetal,klemanparodi,lesliestewartnakagawa,lukyanchuk98,mcmillan},\cite[Chapter 6]{stewart04}). We will restrict attention to smectic A liquid crystals, for which it is often assumed that $r(\x)$ is constant, with the free-energy density being expressed in terms of $\n$ and $\phi$. For example, the free-energy functional proposed in \cite{klemanparodi} is given by
\be 
\label{kp}I(\n,\phi)=\int_\om \left( W(\n,\nabla \n)+\half {\bf B}(\n-\nabla\phi)\cdot(\n-\nabla\phi)\right)\,d\x,
\ee
where ${\bf B}= B_\perp\1+(B_\parallel-B_\perp)\n\otimes\n$ and $B_\perp, B_\parallel$ are positive material constants. A slightly different form for the free energy following \cite{lukyanchuk98} is studied in \cite{calderer08continuum}. In \cite{weinan} it is argued that a good approximation to \eqref{kp} is given by 
\begin{equation}\label{2.14}
 \int_\Omega  (K_1(\textrm{div}\,\n)^2+B_\parallel(|\nabla \phi|-1)^2)\,d\x,
\end{equation}
together with the constraint
\begin{equation}\label{2.13b}
 {\bf n} = \frac{\nabla \phi}{|\nabla \phi|}
\end{equation}
that rigidly enforces that the director points parallel to the normal.  This form is also arrived at in the later paper of \cite{santangelo05curvature}. Because the second term in \eqref{kp} is convex and coercive in $\n-\nabla\phi$ a routine application of the direct method of the calculus of variations shows that $I(\n,\phi)$ attains a minimum over $\n\in W^{1,2}(\om;S^2), \;\nabla\phi\in L^2(\om;\R^3)$ subject to suitable boundary conditions. However it is not so clear how to prove the existence of a minimizer for   the functional \eqref{2.14} under the constraint \eqref{2.13b}.

In seeking a suitable generalization of these models for smectic A which would allow for planar discontinuities in   $\n$, we found it easier instead to follow the route proposed in a 
 recent paper of Pevnyi, Selinger \& Sluckin \cite{selingersluckin}, who criticize the modelling of smectics in terms of the complex order parameter $\Psi$ and the de Gennes free-energy functional on the two grounds  (i) that $\Psi$ cannot be chosen as a single-valued function for index $\half$ defects, for essentially the same reason as discussed in Section \ref{defects} that such defects are not orientable (ii) that it does not represent the local free-energy density at the length-scale of the smectic layers themselves (so that, for example, it cannot model dislocations in the smectic layers). (A related criticism to (i), that the de Gennes free-energy density is not invariant to the transformation $\Psi\to\Psi^*$, is made in \cite{caldererpalffy}.) To remedy these deficiencies Pevnyi, Selinger \& Sluckin propose a modified free-energy density in which the variables are the director $\n$ and the density variation $\rho$. We consider the modified version of this functional
\be
 I(\vQ,\rho) = \int_\Omega \left(\psi_E(\vQ, \nabla \vQ) +B\left|D^2 \rho +\frac{q^2}{3s} \left(3\vQ + s\1 \right) \rho \right|^2 + f(\rho)\right)\,dx &\nonumber\\ + k \int_{S_\vQ} \left| \vQ_+-\vQ_- \right|^r\, d \mathcal{H}^2,& \label{0.22}
\ee
under the uniaxial constraint $\Q(\x)=s\left(\n(\x)\otimes\n(\x)-\frac{1}{3}\1\right)$ for some unit vector $\n(\x)$, where $s>0$ is constant. In \eqref{0.22} $B>0, k>0, q>0, r \in (0,1)$ are constants and $f(\rho)=\frac{a\rho^2}{2}+\frac{b\rho^3}{3} + \frac{c\rho^4}{4}$ for constants $a,b$ and $c>0$. We suppose that the elastic energy $\psi_E$ is a continuous function of its arguments satisfying\vspace{.1in}
 
 \noindent (H1)\;\;\;$\psi_E(\vQ,{\bf H})$ is convex in ${\bf H}$ for any $\vQ$ of the form $s\left( \vn\otimes \vn - \frac{1}{3}\1\right)$,\vspace{.05in}

\noindent (H2)\;\;\;$\psi_E(\vQ,{\bf H}) \geqslant C|{\bf H}|^p + D$ for constants $C>0$, $D$, and $p>1$.\vspace{.1in}

The free energy \eqref{0.22} reduces to that in \cite{selingersluckin} when $\psi_E(\Q)=\frac{K}{2s^2}|\nabla\Q|^2$ and $k=0$, though it is expressed in terms of uniaxial $\Q$ and $\rho$ rather than $\n$ and $\rho$. We allow subquadratic growth of $\psi_E(\Q,\nabla\Q)$ in $\nabla\Q$ so as to be able to model disclinations as line defects with finite energy. As regards the integral over the jump set $S_\Q$, we have chosen the simple form \eqref{fexQ} so as to be able to prove Theorem \ref{ssexist} below. However it would be interesting to allow the integrand to have a nontrivial dependence  on $\nu$ as in \eqref{represent} so as to better represent the effects of dislocations at a curvature wall.   

We now define a set of admissible pairs $\Q,\rho$ by
\be \nonumber
 \mathcal{A}:= \left\{  \Q \in SBV \left( \Omega , \mathbb{R}^{3\times 3} \right),  \rho \in W^{2,2}\left( \Omega,\mathbb{R} \right) :\right. \hspace{1in}&\\ \left. \vQ = s\left(\n\otimes \vn-\frac{1}{3}\1\right),\,\, |\n|=1,\,\, \Q|_{\partial \Omega} = \overline{\Q}\, \right\}&\label{0.23}
\ee
for some sufficiently smooth $\overline{\vQ}:\mathbb{R}^3 \to M^{3\times 3} $. Here $W^{2,2}\left( \Omega,\mathbb{R} \right)$ denotes the space of $\rho\in W^{1,2}(\om;\R)$ with weak derivative $\nabla \rho\in W^{1,2}(\om,\R^3)$. The boundary condition $\Q|_{\partial \Omega} = \overline{\Q}$ is interpreted in the sense described in Section \ref{order} above. 
We can now state
\begin{thm}[\cite{u14}]
\label{ssexist}
 Under the above hypotheses there exists a pair $(\Q^*,\rho^*)$ that minimizes $I(\vQ,\rho)$ over $(\vQ,\rho)\in \mathcal A$.
\end{thm}
This result is of course just a first step towards determining whether such models can be a useful predictive tool for defect morphologies such as those observed in \cite{michel06,coursaultetal2014}.

\section{Concluding remarks}
In this paper we have outlined a possible programme for studying  models for liquid crystals that allow the director to have 2D defects, as well as line and point defects. These models need to be formulated in terms of a suitable function space, for which a natural candidate is $SBV$. In addition to the possible applications  described in Section \ref{planar}, such models may be relevant (see \cite{bedfordcholesterics}) in the high chirality limit for describing cholesteric fingering (see, for example, \cite{smalyukhetal05}) and structures observed in blue phases (see, for example, \cite{wrightmermin}).  Until more work is done on the predictions of such theories it is difficult to estimate their utility.

In the case of order reconstruction, we saw how the Landau - de Gennes energy, which has a five-dimensional order parameter $\Q$ and is expected to have smooth minimizers,  can provide insight into the internal structure of an interface at which, according to the   picture given  by a director theory, the two-dimensional order parameter $\n$ has a discontinuity. This does not imply that director theories involving jumps in $\n$ are not useful. For example, the theory of inviscid gas dynamics is an important model for studying shock waves, precisely because it represents them as singularities, even though introducing viscosity smooths out the shock waves and gives them internal structure.  

That the function space is part of the model seems unquestionable. In this paper we have taken the phenomenological approach of choosing what seems to be the smallest function space capable of describing observed defects. However, it is also natural to seek a deeper physical explanation by obtaining both the energy functional and the function space from passing to a suitable asymptotic limit in a more detailed theory, such as going from Landau - de Gennes to a director model in the limit of small elastic coefficients, or from a molecular model to an Onsager or Landau - de Gennes model in a thermodynamic limit. But this just transfers the question  to that of  where the appropriate function space for the more detailed model comes from. For example, we have seen no evidence of a physical situation in which the $\Q$ tensor has a spatial singularity, but such situations may exist and it does not seem obvious from a theoretical point of view why they should not.

 \section*{Acknowledgements} The research of JMB was supported by 
EPSRC
(GRlJ03466, the Science and Innovation award to the Oxford Centre for Nonlinear
PDE EP/E035027/1, and EP/J014494/1), the European Research Council under the European Union's Seventh Framework Programme
(FP7/2007-2013) / ERC grant agreement no 291053 and
 by a Royal Society Wolfson Research Merit Award.  The research of SJB is supported by an EPSRC CASE studentship with Hewlett-Packard Limited. 

 We are grateful to Fabrice Bethuel,  Christopher Newton, Michaela Nieuwenhuis, Matteo Ricci, Jonathan Selinger, Tim Sluckin, Mark Wilkinson, Arghir Zarnescu    and especially Andrea Braides, Adriana Garroni, Emmanuelle Lacaze, and Claudio Zannoni for helpful suggestions and discussions. The research was initiated during the Isaac Newton Institute research programme on {\it The Mathematics of Liquid Crystals} from January-July 2013, whose support is gratefully acknowledged. The paper was partly written while JMB was visiting the School of Mathematics and Physics, University of Queensland.
 
\bibliography{gen2,balljourn,ballconfproc,ballprep}
\bibliographystyle{abbrv}

\end{document}